\documentclass[12pt]{birkjour}
\usepackage{amsmath}
\usepackage{version}
\usepackage{amsfonts}
\numberwithin{equation}{section}
\newtheorem{theorem}{Theorem}[section]

\newtheorem{example}[theorem]{Example}

\newtheorem{proposition}[theorem]{Proposition}

\newtheorem{remark}[theorem]{Remark}

\newcommand\beq{\begin{equation}}
\newcommand\ds{\displaystyle}
\newcommand\de{\delta}
\newcommand\eeq{\end{equation}}

\newcommand\al{\alpha}

\newcommand\la{\lambda}

\newcommand\si{\sigma}

\newcommand\D{\mathbb D}
\newcommand\T{\mathbb T}
\newcommand\C{\mathbb C}
\newcommand\E{\mathbb E}
\newcommand\G{\mathbb{G}}
\newcommand\R{\mathbb R}

\newcommand\half{{\tfrac{1}{2}}}

\newcommand\df{\stackrel{\rm def}{=}}

\newcommand\nn{\nonumber}

\DeclareMathOperator\diag{diag}

\DeclareMathOperator\tr{tr}

\let\phi\varphi
\begin{document}
\title[Analysable instances of $\mu$-synthesis]{Some analysable instances of mu-synthesis}

\author{ N. J.  Young}
\address{%
Department of Pure Mathematics\\
Leeds University\\
Leeds LS2 9JT\\
England}

\date{17th March 2011}
\dedicatory{To Bill Helton, inspiring mathematician and friend}
\begin{abstract}
I describe a verifiable criterion for the solvability of the $2\times 2$ spectral Nevanlinna-Pick problem with two interpolation points, and likewise for three other special cases of the $\mu$-synthesis problem.
  The problem is to construct an analytic $2\times 2$ matrix function $F$ on the unit disc subject to a finite number of interpolation constraints and a bound on the cost function  $\sup_{\la\in\D} \mu(F(\la))$, where $\mu$ is an instance of the structured singular value.  
\end{abstract}


\maketitle

\section{Introduction} \label {intro}

It is a pleasure to be able to speak at a meeting in San Diego in honour of Bill Helton,  through whose early papers, especially \cite{Hel78}, I first became interested in applications of operator theory to engineering.  I shall discuss a problem of Heltonian character: a hard problem in pure analysis, with immediate applications in control engineering, which can be addressed by operator-theoretic methods.  Furthermore, the main advances I shall describe are based on some highly original ideas of Jim Agler, so that San Diego is the ideal place for my talk.

The $\mu$-synthesis problem is an interpolation problem for analytic matrix functions, a generalization of the classical problems of Nevanlinna-Pick, Carath\'eodory-Fej\'er  and Nehari.  The symbol $\mu$ denotes a type of cost function that generalizes the operator and $H^\infty$ norms, and the $\mu$-synthesis problem is to construct an analytic matrix function $F$ on the unit disc satisfying a finite number of interpolation conditions and such that $\mu(F(\la)) \leq 1$ for $|\la| < 1$.  The precise definition of $\mu$ is in Section \ref{ssv} below, but for most of the paper we need only a familiar special case of $\mu$ -- the spectral radius of a square matrix $A$, which we denote by $r(A)$.

The purpose of this lecture is to present some cases of the $\mu$-synthesis problem that  are amenable to analysis.  I shall summarize some results that are scattered through a number of papers, mainly by Jim Agler and me but also several others of my collaborators, without attempting to survey all the literature on the topic.  I shall also say a little about recent results of some specialists in several complex variables which bear on the matter and may lead to progress on other instances of $\mu$-synthesis.

Although the cases to be described here are too special to have significant practical applications, they do throw some light on the $\mu$-synthesis problem.  More concretely, the results below could be used to provide test data for existing numerical methods and to illuminate the phenomenon (known to engineers) of the numerical instability of some $\mu$-synthesis problems.  

We are interested in citeria for $\mu$-synthesis problems to be solvable.  Here is an example.  We denote by $\D$ and $ \T$ the open unit disc and the unit circle respectively in the complex plane $\C$.
\begin{theorem} \label{main}
Let $\la_1, \la_2\in\D$ be distinct points, let $W_1,W_2$ be nonscalar $2\times 2$ matrices of spectral radius less than $1$ and let $s_j=\tr W_j, \, p_j=\det W_j$ for $j=1,2$.   The following three statements are equivalent:
\begin{enumerate}
\item[(1)]
 there exists an analytic function $F:\D\to\C^{2\times 2}$ such that
\[
F(\la_1)= W_1, \qquad F(\la_2) = W_2
\]
and
\[
r(F(\la)) \leq 1 \quad \mbox{ for all } \la \in \D;
\]
\item[(2)]
\[
\max_{\omega\in\T} \left| \frac{(s_2p_1-s_1p_2)\omega^2 +2(p_2-p_1)\omega + s_1-s_2}{(s_1-\bar s_2p_1)\omega^2 - 2(1-p_1\bar p_2)\omega +\bar s_2-s_1\bar p_2}\right| \leq \left|\frac{\la_1 -\la_2}{1-\bar\la_2\la_1} \right|;
\]
\item[(3)]
\[
\left[ \frac{\overline{(2-\omega s_i)} (2-\omega s_j) - \overline{(2\omega p_i- s_i)}(2\omega p_j-s_j)}{1-\bar \la_i \la_j}\right]_{i,j=1}^2 \geq 0
\]
for all $\omega\in\T$.
\end{enumerate}
\end{theorem}

The paper is organised as follows.   Section \ref{specNP} contains the definition of the spectral Nevanlinna-Pick problem, sketches the ideas that led to Theorem \ref{main} -- reduction to the complex geommetry of the symmetrized bidisc $\G$, the associated ``magic functions" $\Phi_\omega$ and the calculation of the Carath\'eodory distance on $\G$ -- and fills in the final details of the proof of Theorem \ref{main} using the results of \cite{AY5}.  It also discusses ill-conditioning and the possibility of generalization of Theorem \ref{main}.
  In Section \ref{specCF} there is an analogous solvability criterion for a variant of the spectral Nevanlinna-Pick problem in which the two interpolation points coalesce (Theorem \ref{scf}).  In Section \ref{ssv}, besides the definition of $\mu$ and $\mu$-synthesis, there is some motivation and history.  Important work by 
H. Bercovici, C. Foia\c{s} and A. Tannenbaum is briefly described, as is Bill Helton's alternative approach to robust stabilization problems.  In Section \ref{nextcase} we consider an instance of $\mu$-synthesis other than the spectral radius.  Here we can only obtain a solvability criterion in two very special circumstances (Theorems \ref{delE} and \ref{lastcase}).  The paper concludes with some speculations in Section \ref{conclusion}.

We shall denote the  closed unit disc in the complex plane by  $\Delta$.

\section{The spectral Nevanlinna-Pick problem} \label{specNP}
A particularly appealing special case of the $\mu$-synthesis problem is the {\em spectral Nevanlinna-Pick problem}:\\

\noindent {\bf Problem SNP}  {\em  Given distinct points $\la_1,\dots,\la_n \in \D$ and $k\times k$ matrices $W_1,\dots, W_n$, construct an analytic $k\times k$ matrix function $F$ on $\D$ such that
\beq \label{interp}
F(\la_j) = W_j \quad \mbox{ for } j=1,\dots, n
\eeq
and
\beq\label{rleq1}
r(F(\la)) \leq 1\quad \mbox{ for all } \la \in \D.
\eeq }

When $k=1$ this is just the classical Nevanlinna-Pick problem, and it is well known that a suitable $F$ exists if and only if a certain $n \times n$ matrix formed from the $\la_j$ and $W_j$  is positive (this is {\em Pick's Theorem}).  We should very much like to have a similarly elegant solvability criterion for the case that $k > 1$, but strenuous efforts by numerous mathematicians over three decades have failed to find one.

About 15 years ago Jim Agler and I devised a new approach to the problem in the case $k=2$ based on operator theory and a dash of several complex variables (\cite{AY99} to \cite{AY07}).   
Since interpolation of the eigenvalues fails, how about interpolation of the coefficients of the characteristic polynomials of the $W_j$, or in other words of the elementary symmetric functions of the eigenvalues?  This thought brought us to the study of the complex geometry of a certain set $\Gamma \subset \C^2$, defined below.
By this route  we were able to analyse quite fully the simplest then-unsolved case of the spectral Nevanlinna-Pick problem: the case $n=k=2$.  For the purpose of engineering application this is a modest achievement, but it nevertheless constituted progress. It had the merit of revealing some unsuspected intricacies of the problem, and may yet lead to further discoveries.

\subsection{The symmetrized bidisc $\Gamma$}
We introduce the notation
\begin{align}\label{notation}
\Gamma &= \{ (z+w, zw): z,w \in\Delta \}, \\
\G &=\{(z+w,zw):z,w\in\D \}. \nn
\end{align}
$\Gamma$ and $\G$ are called the {\em closed} and {\em open symmetrized bidiscs} respectively.  Their importance lies in their relation to the sets
\begin{align*}
\Sigma  \df \{ A\in \C^{2\times 2} : r(A) \leq 1\}, \\
\Sigma^o  \df \{ A\in \C^{2\times 2} : r(A) < 1\}.
\end{align*}
$\Sigma$ and its interior $\Sigma^o$ are sometines called  ``spectral unit balls", though the terminology is misleading since they are not remotely ball-like, being unbounded and non-convex.  Observe that, for a $2\times 2$ matrix $A$, 
\begin{align*}
A \in \Sigma & \Leftrightarrow \mbox{ the zeros of the polynomial }\la^2 - \tr A \la +\det A \mbox{ lie in } \Delta \\
   & \Leftrightarrow \tr A=z+w,  \, \det A= zw \mbox{ for some }z, \, w \in\Delta.
\end{align*}
  We thus have the following simple assertion. 
\begin{proposition} \label{SigGam}
For any $A\in\C^{2\times 2}$
\begin{align*}
A \in\Sigma &\mbox{  if and only if  } (\tr A, \det A) \in \Gamma,\\
A \in\Sigma^o &\mbox{  if and only if  } (\tr A, \det A) \in \G.
\end{align*}
\end{proposition}
Consequently, if $F:\D \to \Sigma$ is analytic and satisfies the equations (\ref{interp}) above, where $k=2$, then $h\df(\tr F, \det F)$ is an analytic map from $ \D$ to $\Gamma$ satisfying the interpolation conditions
\beq \label{GamInterp}
h(\la_j) = (\tr W_j, \det W_j)\mbox{  for  } j=1,\dots, n.
\eeq
Let us assume that none of the target matrices $W_j$ is a scalar multiple of the identity.  On this hypothesis it is simple to show the converse \cite{BFT1} by similarity transformation of the $W_j$ to companion form.
\begin{proposition}\label{SigtoGam}
Let $\la_1,\dots,\la_n$ be distinct points in $\D$ and let $W_1,\dots,W_n$ be nonscalar $2\times 2 $ matrices.  There exists an analytic map $F:\D\to C^{2\times 2}$ such that equations \eqref{interp} and \eqref{rleq1} hold if and only if there exists an analytic map $h:\D\to \Gamma$ that satisfies the conditions \eqref{GamInterp}.
\end{proposition}
  We have therefore (in the case $k=2$)  reduced the given analytic interpolation problem for $\Sigma$-valued functions to one for $\Gamma$-valued functions (the assumption on the $W_j$ is harmless, since any constraint for which $W_j$ is scalar may be removed by the standard process of Schur reduction).

Why is it an advance to replace $\Sigma$ by $\Gamma$?  For one thing, of the two sets, the geometry of $\Gamma$ is considerably the less rebarbative. $\Sigma$ is an unbounded, non-smooth $4$-complex-dimensional set with spikes shooting off to infinity in many directions.  $\Gamma$ is somewhat better: it is compact and only $2$-complex-dimensional, though $\Gamma$ too is non-convex and not smoothly bounded.  But the true reason that $\Gamma$ is amenable to analysis is that there is a $1$-parameter family of linear fractional functions, analytic on $\G$, that has special properties {\em vis-\`a-vis} $\Gamma$.  For $\omega$ in the unit circle $\T$ we define
\beq \label{defPhi}
\Phi_\omega(s,p) = \frac{2\omega p - s}{2-\omega s}.
\eeq
We use the variables $s$ and $p$ to suggest ``sum" and ``product". 
The $\Phi_\omega$ determine $\G$ in the following sense.
\begin{proposition}\label{critG}
For every $\omega\in\T$, $\Phi_\omega$ maps $\G$
analytically into $\D$.  Conversely, if $(s,p) \in\C^2$ is such that  $|\Phi_\omega(s,p)| < 1$ for all $\omega\in\T$, then $(s,p) \in\G$.
\end{proposition}
Both statements can be derived from the identity
\[
|2-z-w|^2- |2zw-z-w|^2=2(1-|z|^2)|1-w|^2+ 2(1-|w|^2)|1-z|^2.
\]
See \cite[Theorem 2.1]{AY5} for details. 

 There is an analogous statement for $\Gamma$, but there are some subtleties.  For one thing $\Phi_\omega$ is undefined at $(2\bar\omega, \bar\omega^2) \in\Gamma$ when $\omega\in\T$.
\begin{proposition}
For every $\omega\in\T$, $\Phi_\omega$ maps $\Gamma \setminus \{(2\bar\omega, \bar\omega^2)\}$ 
analytically into $\Delta$.  Conversely, if $(s,p) \in\C^2$ is such that  $|\Phi_\omega(rs,r^2p)| < 1$ for all $\omega\in\T$ and $0<r<1$ then $(s,p) \in\Gamma$.
\end{proposition}
In the second statement of the proposition the parameter $r$ is needed: it does not suffice that $|\Phi_\omega(s,p)| \leq 1$ for all $\omega\in\T$ (in the case that $p=1$ the last statement is true if and only if $s\in\R$, whereas for $(s,p)\in\Gamma$, of course $|s| \leq 2$).

  We found the functions $\Phi_\omega$ by applying Agler's theory of families of operator tuples \cite{AY99,AY00}.  We studied the family $\mathcal{F}$ of commuting pairs of operators for which $\Gamma$ is a spectral set, and its dual cone $\mathcal{F}^\perp$ (that is, the collection of hereditary polynomials that are positive on $\mathcal{F}$).  Agler had previously done the analogous analysis for the bidisc, and shown that the dual cone was generated by just two hereditary polynomials; this led to his celebrated realization theorem for bounded analytic functions on the bidisc.  On incorporating symmetry into the analysis we found that the cone $\mathcal{F}^\perp$ had the $1$-parameter family of generators $1-\Phi_\omega^\vee \Phi_\omega, \, \omega\in\T$.  From this fact many conclusions follow: see \cite{AY07} for more on these ideas.

 Operator theory played an essential role in our discovery of the functions $\Phi_\omega$.  Once they are known, however, the geometry of $\G$ and $\Gamma$ can be developed without the use of operator theory.

\subsection{A necessary condition}
Suppose that $F$ is a solution of the spectral Nevanlinna-Pick problem \eqref{interp}, \eqref{rleq1} with $k=2$.  Let us write $s_j=\tr W_j, \, p_j=\det W_j$ for $j=1,\dots,n$.  For any $\omega \in\T$ and $0<t<1$ the composition
\[
\D \stackrel{tF}{\longrightarrow} \Sigma^o \stackrel{(\tr,\det)}{\longrightarrow} \G \stackrel{\Phi_\omega}{\longrightarrow} \D
\]
is an analytic self-map of $\D$ under which
\[
\la_j \mapsto \Phi_\omega( ts_j, t^2p_j)  = \frac{2\omega t^2p_j - ts_j}{2-\omega ts_j}
\quad \mbox{ for } j=1,\dots, n.
\]
Thus, by Pick's Theorem, 
\beq\label{picko}
\left[ \frac{1 - \overline{\Phi_\omega}(ts_i, t^2p_i) \Phi_\omega(ts_j, t^2p_j)}{1-\bar \la_i \la_j} \right]_{i,j=1}^n \geq 0.
\eeq
On conjugating this matrix inequality by $\diag \{ 2-\omega ts_j \}$ and letting $\al=t\omega$ we obtain the following necessary condition for the solvability of a $2\times 2$ spectral Nevanlinna-Pick condition \cite[Theorem 5.2]{AY99}.
\begin{theorem} \label{neccond}
If there exists an analytic map $F: \D\to \Sigma$ satisfying the equations
\[
F(\la_j) = W_j \quad \mbox{ for } j=1,\dots, n
\]
and
\[
r(F(\la)) \leq 1\quad \mbox{ for all } \la \in \D
\]
then, for every $\al$ such that $|\al| \leq 1$,
\beq \label{poscond}
\left[ \frac{\overline{(2-\al s_i)} (2-\al s_j) - |\al|^2\overline{(2\al p_i- s_i)}(2\al p_j-s_j)}{1-\bar \la_i \la_j}\right]_{i,j=1}^n \geq 0
\eeq
where
\[
s_j = \tr W_j, \qquad p_j= \det W_j \quad \mbox{ for } j= 1,\dots, n.
\]
\end{theorem}
In the case that the $W_j$ all have spectral radius strictly less than one, the condition \eqref{poscond} holds for all $\al \in\Delta$ if and only if it holds for all $\al\in\T$, and hence the condition only needs to be checked for a one-parameter pencil of matrices.  It is of course
 less simple than the classical Pick condition in that it comprises an infinite collection of algebraic inequalities, but it is nevertheless checkable in practice with the aid of standard numerical packages.  Its major drawback is that it is {\em not} sufficient for solvability of the $2\times 2$ spectral Nevanlinna-Pick problem.
\begin{example} \label{3ptex} \rm
Let $0<r<1$ and let 
\[
h(\la)=\left( 2(1-r)\frac{\la^2}{1+r\la^3}, \frac{\la(\la^3+r)}{1+r\la^3}\right).
\]
Let $\la_1,\la_2,\la_3$ be any three distinct points in $\D$ and let $h(\la_j) = (s_j,p_j)$ for $j= 1,2, 3$.  We can prove \cite{ALYgamInner} that, in any neighbourhood of $(s_1,s_2,s_3)$ in $(2\D)^3$, there exists a point $(s'_1, s'_2, s'_3)$ such that $(s'_j,p_j)\in \G$, the Nevanlinna-Pick data
\[
\la_j \mapsto \Phi_\omega(s'_j,p_j), \quad j=1,2,3,
\]
are solvable for all $\omega \in \T$, but the Nevanlinna-Pick data
\[
\la_j \mapsto \Phi_{m(\la_j)}(s'_j,p_j), \quad j=1,2,3,
\]
are unsolvable for some Blaschke factor $m$.  It follows that the interpolation data
\[
\la_j \mapsto (s'_j, p_j), \quad j=1,2,3,
\]
satisfy the necessary condition of Theorem \ref{neccond} for solvability, and yet there is no analytic function $h:\D\to\Gamma$ such that $h(\la_j)=(s'_j,p_j)$ for $j=1,2,3$.

\rm Hence, if we choose nonscalar $2\times 2$  matrices $W_1,W_2,W_3$ such that $(\tr W_j,\det W_j)=(s_j,p_j)$, then the spectral Nevanlinna-Pick problem with data $\la_j \mapsto W_j$ satisfies the necessary condition of Theorem \ref{neccond} and yet has no solution.
\end{example}
See also \cite{Bha07} for another example.
\subsection{Two points and two-by-two matrices}
When $n=k=2$ the condition in Theorem \ref{neccond} {\em is} sufficient for the solvability of the spectral Nevanlinna-Pick problem.

We shall now prove the 
 main theorem from Section \ref{intro}.  Recall the statement:\\
{\bf Theorem 1.1.}  
\begin{em}
Let $\la_1, \la_2\in\D$ be distinct points, let $W_1,W_2$ be nonscalar $2\times 2$ matrices of spectral radius less than $1$ and let $s_j=\tr W_j, \, p_j=\det W_j$ for $j=1,2$.   The following three statements are equivalent:
\begin{enumerate}
\item[(1)]
 there exists an analytic function $F:\D\to\C^{2\times 2}$ such that
\[
F(\la_1)= W_1, \qquad F(\la_2) = W_2
\]
and
\[
r(F(\la)) \leq 1 \quad \mbox{ for all } \la \in \D;
\]
\item[(2)]
\beq\label{deleqd}
\max_{\omega\in\T} \left| \frac{(s_2p_1-s_1p_2)\omega^2 +2(p_2-p_1)\omega + s_1-s_2}{(s_1-\bar s_2p_1)\omega^2 - 2(1-p_1\bar p_2)\omega +\bar s_2-s_1\bar p_2}\right| \leq \left|\frac{\la_1 -\la_2}{1-\bar\la_2\la_1} \right|;
\eeq
\item[(3)]
\beq \label{poscond2}
\left[ \frac{\overline{(2-\omega s_i)} (2-\omega s_j) - \overline{(2\omega p_i- s_i)}(2\omega p_j-s_j)}{1-\bar \la_i \la_j}\right]_{i,j=1}^2 \geq 0
\eeq
for all $\omega\in\T$.
\end{enumerate}
\end{em}
  The proof depends on some elementary notions from the theory of invariant distances.  A good source for the general theory is \cite{JaPf}, but here we only need the following rudiments.

We denote by $d$ the pseudohyperbolic distance on the unit disc $\D$:
\[
d(\la_1,\la_2) = \left|\frac{\la_1 -\la_2}{1-\bar\la_2\la_1} \right| \quad\mbox{ for } \la_1,\la_2 \in\D.
\]
 For any domain $\Omega\in\C^n$ we define the {\em Lempert function} $\de_\Omega:\Omega\times \Omega \to \R^+$ by
\beq\label{defLemp}
\de_\Omega(z_1,z_2) = \inf d(\la_1,\la_2)
\eeq
over all $\la_1,\la_2\in\D$ such that there exists an analytic map $h:\D \to \Omega$ such that $h(\la_1)=z_1$ and $ h(\la_2)=z_2$.  We define\footnote{Conventionally the definition of the Carath\'eodory distance contains a $\tanh^{-1}$ on the right hand side of \eqref{defCara}.  For present purposes it is convenient to omit the $\tanh^{-1}$.}  the {\em Carath\'eodory distance} $C_\Omega :\Omega \times \Omega \to \R^+$ by
\beq\label{defCara}
C_\Omega(z_1,z_2) = \sup d(f(z_1),f(z_2))
\eeq
over all analytic maps $f: \Omega \to \D$.  If $\Omega$ is bounded then $C_\Omega$ is a metric on $\Omega$.

It is not hard to see (by the Schwarz-Pick Lemma) that $C_\Omega \leq \de_\Omega$ for any domain $\Omega$.  The two quantities $C_\Omega, \, \de_\Omega$ are not always equal -- the punctured disc provides an example of inequality.  The question of determining the domains $\Omega$ for which $C_\Omega = \de_\Omega$ is one of the concerns of invariant distance theory.
\begin{proof}
Let $z_j =(s_j,p_j) \in\G$.\\
(1)$\Leftrightarrow$(2)
In view of Proposition \ref{SigtoGam} we must show that the inequality \eqref{deleqd} is equivalent to the existence of an analytic $h:\D\to\Gamma$ such that $h(\la_j)=z_j$ for $j=1,2$.  By definition of the Lempert function $\de_\G$, such an $h$ exists if and only if
\[
\de_\G(z_1,z_2) \leq  d(z_1,z_2).
\]
By \cite[Corollary 5.7]{AY5} we have $\de_\G=C_\G$, and by \cite[Theorem 1.1 and Corollary 3.4]{AY5},
\begin{align} \label{CPhi}
C_\G(z_1,z_2) &= \max_{\omega\in\T} d(\Phi_\omega(z_1), \Phi_\omega(z_2))\\
               &= \max_{\omega\in\T} \left| \frac{(s_2p_1-s_1p_2)\omega^2 +2(p_2-p_1)\omega + s_1-s_2}{(s_1-\bar s_2p_1)\omega^2 - 2(1-p_1\bar p_2)\omega +\bar s_2-s_1\bar p_2}\right|. \nn
\end{align}
Thus the desired function $h$ exists if and only if the inequality \eqref{deleqd} holds.\\

\noindent (2)$\Leftrightarrow$(3)   By equation \eqref{CPhi}, the inequality \eqref{deleqd} is equivalent to 
\[
d(\Phi_\omega(z_1), \Phi_\omega(z_2)) \leq d(\la_1,\la_2) \quad\mbox{ for all } \omega\in\T.
\]
By the Schwarz-Pick Lemma, this inequality holds if and only if, for all $\omega\in\T$, there exists a function $f_\omega$ in the Schur class such that $f_\omega(\la_j)=\Phi_\omega(z_j)$ for $j=1,2$.  By Pick's Theorem this in turn is equivalent to the relation
\[
\left[ \frac{1 - \bar\Phi_\omega(z_i) \Phi_\omega(z_j)}{1-\bar \la_i \la_j} \right]_{i,j=1}^2 \geq 0.
\]
Conjugate by $\diag\{2-\omega s_1, 2-\omega s_2\}$ to obtain (2)$\Leftrightarrow$(3).
\end{proof}
\begin{remark} \rm
If one removes the hypothesis that $W_1,W_2$ be nonscalar from Theorem \ref{main} one can still give a solvability criterion.  If both of the $W_j$ are scalar matrices then the problem reduces to a scalar Nevanlinna-Pick problem.  If $W_1 =cI$ and $W_2$ is nonscalar then the corresponding spectral Nevanlinna-Pick problem is solvable if and only if
\[
r((W_2-cI)(I-\bar c W_2)^{-1}) \leq d(\la_1,\la_2)
\]
(see \cite[Theorem 2.4]{AY1}).  This inequality can also be expressed as a somewhat cumbersome algebraic inequality in $c,s_2,p_2$ and $d(\la_1,\la_2)$ \cite[Theorem 2.5(2)]{AY1}.
\end{remark}
\subsection{Ill-conditioned problems}
The results of the preceding subsection suggest that solvability of spectral Nevanlinna-Pick problems depends on the derogatory structure of the target matrices -- that is, in the case of $2\times 2$ matrices, on whether or not they are scalar matrices.  It is indeed so, and in consequence problems in which a target matrix is close to scalar can be very ill-conditioned.
\begin{example} \rm \cite[Example 2.3]{AY1} 
Let $\beta \in\D\setminus \{0\}$ and, for $\al\in\C$ let
\[
W_1(\al) = \left[\begin{array}{cc} 0 & \al \\ 0 & 0 \end{array} \right] , \quad
W_2=\left[\begin{array}{cc} 0 & \beta \\ 0 & \frac{2\beta}{1+\beta} \end{array} \right].
\]
Consider the spectral Nevanlinna-Pick problem with data $0\mapsto W_1(\al), \, \beta \mapsto W_2$.
If $\al=0$ then the problem is not solvable.  If $\al \neq 0$, however, by Proposition \ref{SigtoGam} the problem is solvable if and only if there exists an analytic function $f:\D\to \Gamma$ such that 
\[
f(0)=(0,0) \mbox{ and } f(\beta) = \frac{2\beta}{1+\beta}.
\]
It may be checked \cite{AY2} that
\[
f(\la) = \left(\frac{2(1-\beta)\la}{1-\beta\la}, \frac{\la(\la-\beta)}{1-\beta\la}\right)
\]
is such a function.    Thus the problem has a solution $F_\al$ for any $\al \neq 0$.  Consider a sequence $(\al_n)$ of nonzero complex numbers tending to zero: the functions $F_{\al_n}$ cannot be locally bounded, else they would have a cluster point, which would solve the problem for $\al=0$.  
If $\al$ is, say, $10^{-100}$ then {\em any} numerical method for the spectral Nevanlinna-Pick problem is liable to run into difficulty in this example. 
\end{example}

\subsection{Uniqueness and the construction of interpolating functions}
Problem SNP {\em never} has a unique solution.  If $F$ is a solution of Problem SNP then so is $P^{-1}FP$ for any analytic function $P:\D \to \C^{k\times k}$ such that $P(\la)$ is nonsingular for every $\la\in\D$ and $P(\la_j)$ is a scalar matrix for each interpolation point $\la_j$.  There are always many such $P$ that do not commute with $F$, save in the trivial case that $F$ is scalar.
Nevertheless, the solution of the corresponding interpolation problem for $\Gamma$ {\em can} be unique.  Consider again the case $n=k=2$ with $W_1,W_2$ nonscalar.  By Theorem \ref{main}, the problem is solvable if and only if inequality \eqref{deleqd} holds.  In fact it is solvable {\em uniquely} if and only if inequality \eqref{deleqd} holds with equality.  This amounts to saying that each pair of distinct points of $\G$ lies on a unique complex geodesic of $\G$, which is true by \cite[Theorem 0.3]{AY06}.  (An analytic function $h:\D\to \G$ is a {\em complex geodesic} of $\G$ if  $h$ has an analytic left-inverse).  Moreover, in this case the unique analytic function $h:\D\to \G$ such that $h(\la_j)=(s_j,p_j)$ for $j=1,2$ can be calculated explicitly as follows \cite[Theorem 5.6]{AY5}.

Choose an $\omega_0\in\T$ such that the maximum on the left hand side of \eqref{deleqd} is attained  at $\omega_0$.  Since equality holds in \eqref{deleqd}, we have 
\[
d(\Phi_{\omega_0}(z_1),\Phi_{\omega_0}(z_2)) = d(\la_1,\la_2),
\]
where $z_j=(s_j,p_j)$.  Thus $\Phi_{\omega_0}$ is a Carath\'eodory extremal function for the pair of points $z_1,z_2$ in $\G$.  It is easy (for example, by Schur reduction) to find the unique Blaschke product $p$ of degree at most $2$ such that 
\[
p(\la_1)=p_1, \quad p(\la_2) = p_2\quad \mbox{ and } \quad p(\bar\omega_0) = \bar\omega_0^2).
\]
Define $s$ by
\[
s(\la) = 2\frac{\omega_0 p(\la)-\la}{1-\omega_0 \la} \mbox{  for  } \la\in\D.
\]
Then $h\df(s,p)$ is the required complex geodesic.  

Note that $h$ is a rational function of degree at most $2$.  It can also be expressed in the form of a realization:  $h(\la)= (\tr H(\la), \det H(\la))$ where $H$ is a $2\times 2$ function in the Schur class given by
\[
H(\la) = D+C\la(1-A\la)^{-1} B
\]
for a suitable unitary $3\times 3$ or $4\times 4$ matrix $\left[\begin{array}{cc} A & B\\ C& D\end{array}\right]$ given by explicit formulae (see \cite{AYY}, \cite[Theorem 1.7]{AY06}).

\subsection{More points and bigger matrices}
Our hope in addressing the case $n=k = 2$ of the spectral Nevanlinna-Pick problem was of course that we could progress to the general case.  Alas, we have not so far managed to do so.  We have some hope of giving a good solvability criterion for the case $k=2, \, n=3$, but even the case $n=4$ appears to be too complicated for our present methods.

The case of two points and $k\times k$ matrices, for any $k$, looks at first sight more promising.  There is an obvious way  to generalize the symmetrized bidisc: we define the {\em open symmetrized polydisc }  $\G_k$ to be the domain
\[
\G_k=\{(\si_1(z),\dots, \si_k(z)) : z \in\D^k \} \subset \C^k
\]
where $\si_m$ denotes the elementary symmetric polynomial in $z=(z^1, \dots, z^k)$ for $1\leq m\leq k$.  Similarly one defines the closed symmetrized polydisc $\Gamma_k$.  As in the case $k=2$, one can reduce Problem SNP to an interpolation problem for functions from $\D$ to $\Gamma_k$ under mild hypotheses on the target matrices $W_j$ (specifically, that they be nonderogatory).   However,  the connection between Problem SNP and the corresponding interpolation problems for $\Gamma_k$ are more complicated for $k>2$, because there are more possibilities for the rational canonical forms of the target matrices \cite{NiPfTh}.
The analogues for $\Gamma_k$ of the $\Phi_\omega$ were described by D. J. Ogle \cite{ogle} and subsequently other authors, e.g. \cite{costara,EZ}.  Ogle generalized to higher dimensions the operator-theoretic method of \cite{AY00} and thereby obtained a necessary condition for solvability analogous to Theorem \ref{neccond}.  

The solvability of Problem SNP when $n=2$ is generically equivalent to the inequality
\[
\de_{\G_k}(z_1,z_2) \leq d(\la_1,\la_2)
\]
where $z_j$ is the $k$tuple of coefficients in the characteristic polynomial of $W_j$.  All we need is an effective formula for $\de_{\G_k}$.  It turns out that this is a much harder problem for $k>2$.  In particular, it is {\em false} that $\de_{\G_k} = C_{\G_k}$ when $k>2$.  This discovery \cite{NiPfZw} was disappointing, but not altogether surprising.

There is another type of solvability criterion for the $2\times 2$ spectral Nevanlinna-Pick problem with general $n$ \cite{AY4,Ber03}, but it involves a search over a nonconvex set, and so does not count for the purpose of this paper as an analytic solution of the problem.
Another paper on the topic is \cite{Cos05}.

It is heartening that the study of the complex geometry and analysis of the symmetrized polydisc has been taken up by a number of specialists in several complex variables, including  G. Bharali, C. Costara, A. Edigarian, M. Jarnicki, L. Kosinski, N. Nikolov, P. Pflug, P. Thomas and W. Zwonek.  Between them they have made many interesting discoveries about these and related domains.  There is every hope that some of their results will throw further light on the spectral Nevanlinna-Pick problem.

\section{The spectral Carath\'eodory-Fej\'er problem} \label{specCF}
This is the problem that arises from the spectral Nevanlinna-Pick problem when the interpolation points coalesce at $0$.

\noindent {\bf Problem SCF} {\em Given $k\times k$ matrices $V_0,V_1, \dots, V_n,$ find an analytic function $F: \D\to\C^{k\times k}$ such that
\beq\label{derivs}
F^{(j)}(0) = V_j \quad \mbox{ for } j=0,\dots, n
\eeq
and
\beq\label{specrad}
r(F(\la)) \leq 1 \quad \mbox{ for all } \la\in\D.
\eeq}

This problem also can be converted to an interpolation problem for analytic functions from $\D$ into $\Gamma_k$ \cite[Theorem 2.1]{HMY}, \cite{NiPfTh}.  However, the resulting problem is again hard when $k \geq 2$, and the only truly explicit solution we have is in the case $k=2, n=1$ \cite[Theorem 1.1]{HMY}.
\begin{theorem}\label{scf}
Let
\[
V_m = [ v^m_{ij}]_{i,j=1}^2 \quad \mbox{ for } m= 0, 1
\]
and suppose that $V_0$ is nonscalar.  There exists an analytic function $F:\D \to \C^{2\times 2}$ such that 
\beq\label{scfproblem}
F(0)=V_0, \quad F'(0) = V_1 \quad \mbox{ and } \quad r(F(\la))  < 1 \mbox{ for all } \la\in\D
\eeq
if and only if
\beq\label{carametric}
\max_{|\omega|=1} \left| \frac {(s_1p_0-s_0p_1)\omega^2 +2\omega p_1 - s_1}{\omega^2 (s_0-\bar s_0 p_0) - 2\omega (1-|p_0|^2) + \bar s_0 - s_0 \bar p_0}\right| \leq 1,
\eeq
where 
\begin{align*}
 s_0 &= \tr V_0, \quad p_0= \det V_0, \\
s_1 &= \tr V_1, \quad p_1 = \left| \begin{array}{cc} v^0_{11} & v^1_{12} \\ v^0_{21} & v^1_{22} \end{array} \right| +
   \left| \begin{array}{cc} v^1_{11} & v^0_{12} \\ v^1_{21} & v^0_{22} \end{array} \right|.
\end{align*}
\end{theorem} 
The proof of this theorem in \cite{HMY} again depends on the calculation in \cite{AY5} of the Carath\'eodory metric on $\G$, but this time on the infinitesimal version $c_\G$ of the metric: the left hand side of inequality \eqref{carametric} is the value of $c_\G$ at $(s_0,p_0)$ in the direction $(s_1,p_1)$.  This fact is \cite[Corollary 4.4]{AY5}, but unfortunately there is an $\omega$ missing in the statement of Corollary 4.4.  The proof shows that the correct formula is as in \eqref{carametric}.  An important step is the proof that the infinitesimal Carath\'eodory and Kobayashi metrics on $\G$ coincide.

The ideas behind Theorem \ref{scf} can be used to find solutions of Problem SCF: see \cite[Section 6]{HMY}.  The ideas can also be used to derive a necessary condition for the spectral Carath\'eodory-Fej\'er problem \eqref{derivs}, \eqref{specrad} in the case that $n=1$ and $k>2$ \cite[Theorem 4.1]{HMY}, but there is no reason to expect this condition to be sufficient.

\section{The structured singular value}\label{ssv}

The {\em structured singular value} of a matrix relative to a space of matrices was introduced by J. C. Doyle and G. Stein in the early 1980s \cite{Do82,DoSt} and was denoted by $\mu$.  It is a refinement of the usual operator norm of a matrix and is motivated by the problem of the robust stabilization of a plant that is subject to structured uncertainty.  Initially, in the $H^\infty$ approach to robustness, the uncertainty of a plant was modelled by a meromorphic matrix function (on a disc or half plane) that is subject to an $L^\infty$ bound but is otherwise completely unknown.  The problem of the simultaneous stabilization of the resulting collection of plant models could then be reduced to some classical analysis and operator theory, notably to the far-reaching results of Adamyan, Arov and Krein from the 1970s \cite{Fr87}.

In practice one may have some structural information about the uncertainty in a plant -- for example, that  certain entries are zero.  By incorporating such structural information one should be able to achieve a less conservative stabilizing controller.  The structured singular value was devised for this purpose.  A good account of these notions is in \cite[Chapter 8]{DuPa}.    
Unfortunately, the behaviour of $\mu$ differs radically from that of the operator norm -- for one thing, $\mu$ is not in general a norm at all, and none of the relevant classical theorems (such as Pick's theorem) or methods appear to extend to the corresponding questions for $\mu$.  This provides a challenge for mathematicians: we should help out our colleagues in engineering by creating an AAK-type theory for $\mu$.

For any $A\in \C^{k\times\ell}$ and any subspace $E$ of $\C^{\ell\times k}$ we define the structured singular value $\mu_E(A)$ by
\beq\label{defmu}
\frac{1}{\mu_E(A)} = \inf \{ \|X\|: X \in E,  \, 1-AX \mbox{ is singular} \}
\eeq
with the understanding that $\mu_E(A)=0$ if $1-AX$ is always nonsingular.  

Two instances of the structured singular value are the operator norm $\|.\|$ (relative to the Euclidean norms on $\C^k$ and $ \C^\ell$) and the spectral radius $r$.  If we take $E= \C^{\ell\times k}$ then we find that $\mu_E(A) = \|A\|$.  On the other hand, if $k=\ell$ and we choose $E$ to be the space of scalar multiples of the identity matrix, then $\mu_E(A)=r(A)$.  These two special $\mu$s are in a sense extremal: it is always the case, for any $E$,  that 
 $\mu_E(A) \leq \|A\|$.    If $k=\ell$ and $E$ contains the identity matrix, then $\mu_E(A) \geq r(A)$.  
  A comprehensive discussion of the properties of $\mu$ can be found in \cite{PaDo}.

Here is a formulation of the {\em $\mu$-synthesis problem} \cite{DoSt,DuPa}.

\begin{em}
Given positive integers $k, \ell$, a subspace $E$ of $\C^{\ell\times k}$ and analytic functions $A, B, C$ on $\D$ of types $k\times \ell, k\times k$ and $\ell \times \ell$ respectively, construct an analytic function $F:\D\to \C^{k\times \ell}$ of the form
\beq\label{modelmatch}
F=A+BQC \quad \mbox{ for some analytic }\quad Q:\D\to \C^{k\times \ell}
\eeq
such that
\beq\label{mubound}
\mu_E(F(\la)) \leq 1 \quad \mbox{ for all } \quad \la\in\D.
\eeq
\end{em}

The condition \eqref{modelmatch}, that $F$ be expressible in the form $A+BQC$ for some analytic $Q$,  can be regarded as an interpolation condition on $F$.  In the event that $k=\ell$, $B$ is the scalar polynomial
\[
B(\la) = (\la-\la_1)\dots (\la-\la_n) I
\]
with distinct zeros $\la_j\in\D$ and $C$ is constant and equal to the identity, then  $F$ is expressible in the form $A+BQC$ if and only if
\[
F(\la_1) = A(\la_1), \dots, F(\la_n) = A(\la_n).
\]
With this choice of $B$ and $C$, if we take $E$ to be the space of scalar matrices, we obtain precisely the spectral Nevanlinna-Pick problem.  If we now replace $B$ by the polynomial $\la^n$, we get the spectral Carath\'eodory-Fej\'er problem.

 In engineering applications $\mu$-synthesis problems arise after some analysis is carried out on the plant model to produce the $A, B$ and $C$ in condition \eqref{modelmatch}, and the resulting $B$ and $C$ will not usually be scalar functions.  Nevertheless, explicit pointwise interpolation conditions provide a class of easily-formulated test cases, and it is arguable that such problems are the {\em hardest} cases of $\mu$-synthesis.

Conditions of the form \eqref{modelmatch} are said to be of {\em model matching type} \cite{Fr87}. 

The most sustained attempt to develop an AAK-type theory for the structured singular value in full generality is due to H. Bercovici, C. Foia\c{s} and A. Tannenbaum (\cite{BFKT} to\cite{BFT6}).  They have a far-reaching theory: {\em inter alia} they have constructed many illuminating examples, found properties of extremal solutions and obtained a type of solvability criterion for $\mu$-synthesis problems.  The criterion results from a combination of the Commutant Lifting Theorem with the application of similarity transformations.  To apply the criterion to a concrete spectral Nevanlinna-Pick problem one must solve an optimization problem over a high-dimensional unbounded and non-convex set. We can certainly hope that this is not the last word on the subject of solvability.  Despite the achievements of Bercovici, Foia\c{s} and Tannenbaum, there is still plenty of room for further study of $\mu$-synthesis.

 One of their examples \cite[Section 7, Example 5]{BFT3} exhibits an important fact about the spectral Nevanlinna-Pick problem: {\em diagonalization does not work}. It shows that diagonalization of the target matrices $W_j$ in Problem SNP by similarity transformations, even when possible, does not help solve the problem.  One could hope that if the $W_j$ were diagonal one might be able to  decouple the problem into a series of scalar interpolation problems, but they show that such a hope is vain.

Bill Helton himself, along with collaborators, has developed an alternative approach to the refinement of $H^\infty$ control; his viewpoint is set out in \cite{Hel87}.  His part in the introduction of the results of Adamyan, Arov, Krein and other operator-theorists into robust control theory in the early 1980s is well known.   He subsequently worked extensively (with Orlando Merino, Trent Walker and others) during the 1990s on the more delicate optimization problems that arise from refinements of the basic $H^\infty$ picture of modelling uncertainty.  As in the $\mu$ approach, the aim is to incorporate more subtle specifications and robustness conditions into methods for controller design.
    He developed a very flexible  formulation of such problems as optimization problems over spaces of vector-valued analytic functions on the disc, and devised an algorithm for their numerical solution -- see \cite{HelMerWa93} and several other papers. The authors proved convergence results and described numerical trials.    However, the spectral Nevanlinna-Pick problem cannot be satisfactorily treated by the Helton scheme.  Although it can be cast in the basic problem formulation \cite[Chapter 2]{Hel87}, solution algorithms require smoothness properties (of the function ``$\Gamma$") which the spectral radius does not possess.

\section{The next case of $\mu$}\label{nextcase}
After the two extremes $\mu=\|.\|_{H^\infty}$ and $\mu=r$ the next natural case to consider is the one in which, in \eqref{defmu},  $k=\ell$ and $E$ is the space $\mathrm{Diag}(k)$ of diagonal matrices.  For the rest of this section $\mu$ will denote $\mu_ {\mathrm{Diag}(2)}$ and we shall study the following problem:

{\em Given distinct points $\la_1,\dots,\la_n \in\D$ and $2\times 2$ matrices $W_1,\dots,W_n$, construct an analytic function $F:\D \to \C^{2\times 2}$ such that 
\beq\label{interpnext}
F(\la_j) = W_j \quad \mbox{ for } j= 1,\dots, n
\eeq
and
\beq\label{muleq1}
\mu(F(\la)) \leq 1 \quad \mbox{ for all } \la \in\D.
\eeq}

For the $2\times 2$ spectral Nevanlinna-Pick problem we had some modest success through reduction to an interpolation problem for $\Gamma$-valued functions.  In the present case we tried an analogous approach, with still more modest success \cite{alaa, awy, Y08}.   The following result is \cite[Theorem 9.4 and Remark 9.5(iii)]{awy}.
\begin{theorem} \label{delE}
Let $\la_0\in\D, \, \la\neq 0$, let $\zeta\in\C$ and let
\beq\label{defW}
W_1= \left[ \begin{array}{cc} 0 & \zeta \\ 0 & 0 \end{array}\right], \qquad
W_2 = \left[ \begin{array}{cc} a & * \\ * & b \end{array}\right].
\eeq
Suppose that $|b| \leq |a|$ and let $p=\det W_2$.   There exists an analytic function $F:\D\to \C^{2\times 2}$ such that
\beq\label{muinterp2}
F(\la_1)=W_1, \quad F(\la_2)= W_2 \quad \mbox{ and } \quad \mu(F(\la)) \leq 1 \quad \mbox{ for all } \la \in\D
\eeq
if and only if $|p| < 1$ and 
\[
\left \{ \begin{array}{ll}
 \ds  \frac{|a-\bar b p| + |ab-p|}{1-|p|^2} \leq |\la_0| & \mbox{ if }\zeta \neq 0 \\
      &   \\
 |\la_0|^4 - (|a|^2+|b|^2+2|ab-p|)|\la_0|^2 +|p|^2 \geq 0 \quad & \mbox{ if }\zeta = 0.
 \end{array} \right.
\]
\end{theorem}
The stars in the formula for $W_2$ in \eqref{defW}  denote arbitrary complex numbers.

What is the analog of $\Gamma$ for this case of $\mu$?
To determine whether a $2\times 2$ matrix $A=[a_{ij}]$ satisfies $r(A) \leq 1$ one needs to know only the two numbers $\tr A$ and $\det A$; this fact means that the spectral Nevanlinna-Pick problem can generically be reduced to an interpolation problem for $\Gamma$.  To determine whether $\mu(A)\leq 1$ one needs to know the three numbers $a_{11}, a_{22}, \det A$.  This led us to introduce a domain $\E$ which we call the {\em tetrablock}:
\beq\label{defE}
\E = \{x\in\C^3: 1-x^1 z- x^2 w+ x^3 zw \neq 0 \mbox{ whenever } |z| \leq 1, |w| \leq 1 \}.
\eeq
Its closure is denoted by $\bar \E$.  The name reflects the fact that the intersection of $\E$ with $\R^3$ is a regular tetrahedron.  The domain $\E$ is relevant because $\mu(A) < 1$ if and only if $(a_{11}, a_{22}, \det A) \in \E$.     There exists a solution of the $2$-point $\mu$-synthesis problem \eqref{muinterp2} if and only if the corresponding interpolation problem for analytic functions from $\D$ to $\E$ is solvable \cite[Theorem 9.2]{awy}, and accordingly the solvability problem for this $\mu$-synthesis problem is equivalent to the calculation of the Lempert function $\de_\E$.  As far as I know no one has yet computed $\de_\E$ for a general pair of points of $\E$, but we did calculate it in the case that one of the points is the origin in $\C^3$, that is, we proved a Schwarz lemma for $\E$.   The result is Theorem \ref{delE}.

Observe that ill-conditioning appears in this instance of $\mu$-synthesis too \cite[Remark 9.5(iv)]{awy}. If, in Theorem \ref{delE}, $a=b=p=\half$ then there exists a solution $F_\zeta$ of the problem if and only if
\[
|\la_0| \geq \left\{ \begin{array}{cl} \tfrac 23 & \mbox{ if } \zeta\neq 0 \\
   &  \\
   \frac{1}{\sqrt{2}} & \mbox{ if } \zeta =0
 \end{array} \right.
\]
Thus if $\tfrac 23 < |\la_0| < \frac{1}{\sqrt{2}}$, the $F_\zeta$ are not locally bounded as $\zeta \to 0$, and so are sensitive to small changes in $\zeta$ near $0$.

The complex geometry of $\E$ has also proved to be of interest to researchers in several complex variables.  To my surprise, it  was recently shown \cite{EdKoZw10} that the Lempert function and the Carath\'eodory distance  on $\E$ coincide.  This might be a step on the way to the derivation of a formula for $\de_\E$.  It would suffice to compute $\de_\E$ in the case that one of the two points is of the form $(0,0,\la)$ for some $\la\in[0,1)$, since every point of $\E$ is the image of such a point under an automorphism of $\E$ \cite[Theorem 5.2]{Y08}.

The fourth and final special case of $\mu$-synthesis in this paper is the $\mu$-analog of the $2\times 2$ Carath\'eodory-Fej\'er problem:

{\em Given $2\times 2$ matrices $V_0, \dots, V_n$, construct an analytic function $F:\D\to \C^{2\times 2}$ such that 
\[
F^{(j)}(0) = V_j \mbox{ for } j=0, \dots, n \quad \mbox{ and } \quad \mu(F(\la)) \leq 1 \quad \mbox{ for all } \la\in\D.
\]   }

Again the problem can be reduced to an interpolation problem for $\E$, but the resulting problem has only been solved in an exceedingly special case.
\begin{theorem}\label{lastcase}
Let $V_0,V_1$ be $2\times 2$ matrices such that 
\[
V_0=\left[\begin{array}{cc} 0 & \zeta \\ 0 & 0 \end{array} \right]
\] 
 for some $\zeta\in\C$ and $V_1=  [v_{ij}]$ is nondiagonal.  There exists an analytic function $F:\D\to \C^{2\times 2}$ such that
\[
F(0)=V_0, \quad F'(0)=V_1 \quad \mbox{ and } \quad \mu(F(\la)) \leq 1 \mbox{ for all } \la\in\D
\]
if and only if
\[
\max\{ |v_{11}|, |v_{22}|\} + | \zeta v_{21} | \leq 1.
\]
\end{theorem}
This result follows from \cite[Theorem 2.1]{Y08}.
\section{Conclusion}\label{conclusion}
Although $\mu$-analysis remains a useful tool, it is fair to say that $\mu$-synthesis, as a major technique for robust control system design, has been something of a disappointment up to now.  The trouble is that the $\mu$-synthesis problem is difficult.  It is a highly non-convex problem.  There do exist heuristic numerical methods for addressing particular $\mu$-synthesis problems, notably a Matlab toolbox \cite{matlab} based on the ``$DK$ algorithm" \cite[Section 9.3]{DuPa}, but there is no practical solvability criterion, no fast algorithm nor any convergence theorem for any known algorithm.  For these reasons engineers have largely turned to other approaches to robust stabilization over the past 20 years.  If, however, a satisfactory analytic theory of the problem is developed, engineers' attention may well return to $\mu$-synthesis as a promising design tool.  We are still far from having such a theory, but perhaps these special cases and the interest of the several complex variables community may yet lead to one.

\end{document}